\documentclass[10pt]{amsart}       
 
\title[  Formality  hard Lefschetz property  ]
{ Formality and  hard Lefschetz property  of aspherical manifolds}

\author{Hisashi Kasuya}
\usepackage{amssymb}   
\usepackage{amsmath} 
\usepackage{amscd}
\usepackage{amstext}
\usepackage{amsfonts}
\usepackage[all]{xy}

\address[H.kasuya]{Graduate school of mathematical science university of tokyo japan }
\curraddr{}
\email{khsc@ms.u-tokyo.ac.jp}

\thanks{}

\newcommand{\C}{\mathbb{C}}
\newcommand{\R}{\mathbb{R}} 

\newcommand{\Q}{\mathbb{Q}}
\newcommand{\Z}{\mathbb{Z}}
\newcommand{\g}{\frak{g}}
\newcommand{\n}{\frak{n}}
\newcommand{\U}{{\bf U}}
\newcommand{\Aut}{\rm Aut}

\theoremstyle{remark}
\newtheorem{example}{\rm Example}
\theoremstyle{plain}
\newtheorem{theorem}{Theorem} [section]
\theoremstyle{remark}
\newtheorem{remark}{\rm Remark}
\theoremstyle{plain}
\newtheorem{lemma}[theorem]{Lemma}
\theoremstyle{definition}
\newtheorem{definition}[theorem]{\rm Definition}
\theoremstyle{plain}
\newtheorem{proposition}[theorem]{Proposition}
\theoremstyle{plain}
\newtheorem{corollary}[theorem]{Corollary}
\begin{document} 
\begin{abstract}
For a  Lie group $G=\R^{n}\ltimes _{\phi}\R^{m}$ with the semi-simple action $\phi:\R^{n}\to {\rm Aut}(\R^{m}) $, we show that if $\Gamma$ is  a finite extension of a lattice of $G$
then $K(\Gamma, 1)$ is formal.
Moreover we show that a compact symplectic aspherical manifold with the fundamental group $\Gamma$ satisfies the hard Lefschetz property.
By those results we give many examples of formal solvmanifolds satisfying the hard Lefschetz property but not admitting K\"ahler structures.
 \end{abstract}
\maketitle

\section{Introduction}

Formal spaces(see Definition \ref{df}) in the sense of Sullivan are important in   de Rham homotopy theory.
Well-known examples of formal spaces are compact K$\ddot {\rm a}$hler manifolds (see \cite{DGMS}).
Suppose $\Gamma$ is a torsion-free finitely generated nilpotent group.
Then $K(\Gamma, 1)$ is formal if and only if $\Gamma$ is abelian by Hasegawa's theorem in \cite{H}. 
But in case $\Gamma$ is a virtually polycyclic(see Definition \ref{a-d-1}) group, the formality of $K(\Gamma, 1)$ is more complicated.
One of the purposes of this paper is to apply the way of the algebraic hull of $\Gamma$ to study the formality of $K(\Gamma, 1)$.
For a torsion-free virtually polycyclic group $\Gamma$, we have a unique algebraic group ${\bf H}_{\Gamma}$ with an injective homomorphism $\psi:\Gamma\to {\bf H}_{\Gamma}$ so that:\\
(1)  \ $\psi (\Gamma)$ is Zariski-dense in ${\bf H}_{\Gamma}$.\\
(2) \    The centralizer $Z_{{\bf H}_{\Gamma}}({\bf U}({\bf H}_{\Gamma}))$ of ${\bf U}({\bf H}_{\Gamma})$ is contained in  ${\bf U}({\bf H}_{\Gamma})$.  \\
(3) \ $\dim {\bf U}({\bf H}_{\Gamma})= {\rm rank}\,\Gamma$.   \\
Such ${\bf H}_{\Gamma}$ is called the algebraic hull of $\Gamma$.
We call the unipotent radical of ${\bf H}_{\Gamma}$ the unipotent hull of $\Gamma$ and denote it by ${\bf U}_{\Gamma}$.
In \cite{B}, Baues constructed a compact aspherical manifold $M_{\Gamma}$ with the fundamental group $\Gamma$ which is called the standard $\Gamma$-manifold by the algebraic hull of $\Gamma$.
And he  gave the way of computation of the de Rham cohomology of $M_{\Gamma}$.
By using these results, we prove:
\begin{proposition}\label{abu}
If the  unipotent hull ${\bf U}_{\Gamma}$ of $\Gamma$ is abelian,  $K(\Gamma, 1)$ is formal.
\end{proposition}
So we would like to know criteria for ${\bf U}_{\Gamma}$ to be abelian.
We prove the following theorem.
\begin{theorem}\label{thoo}
Let $\Gamma$ be a torsion-free virtually polycyclic group.
Then the following two conditions are equivalent:\\
$(1)$ ${\bf U}_{\Gamma}$ is abelian.\\
$(2)$ $\Gamma$  is a finite extension group of a lattice of a Lie group  $G=\R^{n}\ltimes_{\phi} \R^{m}$ such that the action $\phi:\R^{n}\to {\rm  Aut}(\R^{m})$ is semi-simple.
\end{theorem}
Therefore we have:
\begin{corollary}\label{mmo} 
If  $\Gamma$ satisfies the condition $(2)$ in Theorem \ref{thoo},
then $K(\Gamma, 1) $ is formal.
\end{corollary}
\begin{remark}
A lattice $\Gamma$ of  $G=\R^{n}\ltimes_{\phi} \R^{m}$ is the form $\Gamma^{\prime}\ltimes_{\phi} \Gamma^{\prime\prime}$ such that $\Gamma^{\prime}$ and $\Gamma^{\prime\prime}$ are lattices of $\R^{n}$ and $\R^{m}$ respectively and the action $\phi$ of $\Gamma^{\prime}$ preserves $\Gamma^{\prime\prime}$.
\end{remark}

As well as formality the hard Lefschetz property(see Definition \ref{dhd})  is an  important property of a compact K$\ddot {\rm a}$hler manifold.
We have the following proposition.
\begin{proposition}\label{LEE}
Let $M$ be a compact symplectic aspherical manifold with the torsion-free virtually polycyclic fundamental group $\Gamma$.
If the unipotent hull ${\bf U}_{\Gamma}$ is abelian,  then $M$ satisfies the hard Lefschetz property.
\end{proposition}
Hence we have:
\begin{corollary}
If  $\Gamma$ satisfies the condition $(2)$ in Theorem \ref{thoo},
then a compact symplectic aspherical manifold with the fundamental gorup $\Gamma$ satisfies the hard Lefschetz property.
\end{corollary}
In \cite{BG}, Benson and Gordon showed that a compact symplectic aspherical manifold with the torsion-free nilpotent fundamental group $\Gamma$ satisfies the hard Lefschetz property if and only if $\Gamma$ is abelian.

As we see in \cite{H} and \cite{BG}, formality and the hard Lefschetz property are strong criteria for aspherical manifolds to admit K\"ahler structures.
But by the results of this paper, we can  obtain many non-K\"ahler formal aspherical manifolds satisfying the hard Lefschetz property.

Let $M$ be a compact aspherical manifold with the  virtually polycyclic fundamental group.
In \cite{BC}, Baues and Cort\'es showed that if $M$ admits a K$\ddot{{\rm a}}$hler structure then the fundamental group of $M$ is virtually abelian(this result is an extension of the result in \cite{A} and \cite{Hn}).
Let $G$ be a simply connected solvable Lie group.
We say that  $G$ is of type (I) if for any $g\in G$  all eigenvalues of the adjoint operator ${\rm Ad}_{g}$ have absolute value $1$.
In \cite{Aus} it was proved that a lattice of a simply connected solvable Lie group $G$ is virtually nilpotent if and only if $G$ is type (I).
Hence we have:
\begin{corollary}
Let  $\Gamma$  be a finite extension group of a lattice of a Lie group  $G=\R^{n}\ltimes_{\phi} \R^{m}$ such that the action $\phi:\R^{n}\to {\rm  Aut}(\R^{m})$ is semi-simple and $G$ is not of type (I).
Then a compact aspherical manifold $M$ with the fundamental group $\Gamma$ is formal but admits no K\"ahler structure.
If $M$ admits a symplectic structure, then $M$ satisfies the hard Lefschetz property.
\end{corollary}
\begin{remark}
In \cite{Hn}, Hasegawa  showed that  a simply connected solvable Lie group $G$ with a virtually abelian lattice such that $G/\Gamma$ admits K\"ahler structure can be written as $G=\R^{2k}\ltimes _{\phi}\C^{l}$ such that
\[\phi(t_{j})((z_{1},\dots ,z_{l}))=(e^{\sqrt{-1}\theta^{j}_{1}t_{j}}z_{1},\dots ,e^{\sqrt{-1}\theta^{j}_{l}t_{j}}z_{l}),
\]
where each $e^{\sqrt{-1}\theta^{j}_{i}}$ is a root of unity.
\end{remark}

 Solvmanifolds are homogeneous spaces of  connected solvable Lie groups.
These are examples of aspherical manifolds with the polycyclic fundamental groups.
In particular for a simply connected solvable Lie group $G$ with a lattice $\Gamma$, the solvmanifold $G/\Gamma$ is a compact aspherical  manifold with the fundamental group $\Gamma$.
As  generalizations of solvmanifolds we define infra-solvmanifolds.
Let $G$ be a simply connected solvable Lie group.
Consider the group ${\rm Aut}(G)\ltimes G$ of affine transformations of $G$ and the projection $p:{\rm Aut}(G)\ltimes G\to {\rm Aut}(G)$.
An infra-solvmanifold is a manifold of the form $G/\Delta$ for a torsion-free subgroup $\Delta$ of ${\rm Aut}(G)\ltimes G$ such that $p(\Delta)$ is contained in a compact subgroup of ${\rm Aut}(G)$.
In \cite{B} Baues showed that every compact infra-solvmanifold is diffeomorphic to a standard $\Gamma$-manifold and for any torsion-free virtually polycyclic group $\Gamma$ the standard $\Gamma $-manifold is diffeomorphic to an infra-solvmanifold $G/\Gamma$ such that $\Gamma\subset {\rm Aut}(G)\ltimes G$ is  a discrete subgroup and $p(\Gamma)$ is finite.
Thus for any $\Gamma$ satisfying the condition $(2)$ in Theorem \ref{thoo} we have a compact infra-solvmanifold $G/\Gamma $ for some $G=\R^{n}\ltimes_{\phi} \R^{m}$ such that the action $\phi:\R^{n}\to {\rm  Aut}(\R^{m})$ is semi-simple.

{\em Notations and teminology}:
Let $k$ be a subfield of $\C$.
A group $\bf G$ is called $k$-algebraic group if $\bf G$ is a Zariski-closed subgroup of $GL_{n}(\C)$ which is defined by polynomials with coefficients in $k$.
Let  ${\bf G}(k)$ denote the set  of  $k$-points of $\bf G$ 
and ${\bf U}({\bf G})$ the maximal Zariski-closed unipotent normal $k$-subgroup of $\bf G$ called the unipotent radical of $\bf G$.
A general reference  is \cite{Bor}.
In this paper, algebraic groups are always written in the bold face.

\section{Algebraic hulls} 
In this section we explain the algebraic hulls of polycyclic groups or simply connected solvable Lie groups.
\begin{definition}\label{a-d-1}
A group $\Gamma$ is {\em polycyclic} if it admits a sequence 
\[\Gamma=\Gamma_{0}\supset \Gamma_{1}\supset \cdot \cdot \cdot \supset \Gamma_{k}=\{ e \}\]
of subgroups such that each $\Gamma_{i}$ is normal in $\Gamma_{i-1}$ and $\Gamma_{i-1}/\Gamma_{i}$ is cyclic.
We set ${\rm rank}\,\Gamma=\sum_{i=1}^{i=k} {\rm rank}\, \Gamma_{i-1}/\Gamma_{i}$ which is  independent of the choice of a sequence $\Gamma_{i}$.
\end{definition}
There are close relations between polycyclic groups and solvable Lie groups.
\begin{theorem}{\rm (\cite[Proposition 3.7, Theorem 4.28]{R})}\label{subpo}
Let $G$ be a simply connected solvable Lie group and $\Gamma$ a lattice in $G$.
Then $\Gamma$ is torsion-free polycyclic and $ \dim G={\rm rank}\,\Gamma$.
Conversely every polycyclic group admits a finite index normal subgroup  which is isomorphic to a lattice in a simply connected solvable Lie group.
\end{theorem}

Let $\Gamma$ be a virtually polycyclic group and $\Gamma^{\prime}$ be a finite index polycyclic subgroup.
We set ${\rm rank}\,\Gamma={\rm rank}\,\Gamma^{\prime}$.

\begin{definition}\label{a-d-2}
Let $k$ be a subfield $\C$.
Let $\Gamma$ be a  torsion-free virtually polycyclic group(resp. simply connected solvable Lie group).
Then a $k$-algebraic group  ${\bf H}_{\Gamma}$ is a {\em $k$-algebraic hull} of  $\Gamma$ if there exists an injective  homomorphism $\psi :\Gamma\to {\bf H}_{\Gamma}(k)$ 
and ${\bf H}_{\Gamma}$ satisfies the following conditions:
\\
(1)  \ $\psi (\Gamma)$ is Zariski-dense in $\bf{H}_{\Gamma}$.\\
(2) \   $Z_{{\bf H}_{\Gamma}}({\bf U}({\bf H}_{\Gamma}))\subset {\bf U}({\bf H}_{\Gamma})$.\\
(3) \ $\dim {\bf U}({\bf H}_{\Gamma})$=${\rm rank}\,\Gamma$ (resp. $\dim \Gamma$).   
\end{definition}
\begin{theorem}\label{a-t-1}{\rm (\cite[Theorem A.1, Corollary A.3]{B})}{\rm (\cite[Proposition 4.40, Lemma 4.41]{R})}
Let $\Gamma$ be a  torsion-free virtually polycyclic group(resp. simply connected solvable Lie group).
Then there exists a $\Q$-algebraic(resp. $\R$-algebraic) hull of $\Gamma$  and for any subfield $k\subset \C$ which contains $\Q$(resp. $\R$) a $k$-algebraic hull of $\Gamma$  is unique up to $k$-algebraic group isomorphism.
\end{theorem}

We call the unipotent radical of ${\bf H}_{\Gamma}$ the unipotent hull of $\Gamma$ and denote it by ${\bf U}_{\Gamma}$.

\begin{lemma}\label{a-l-1} 
Let $\Gamma$ be a torsion-free virtually polycyclic group
and  $\Delta $  a finite index subgroup of $\Gamma$.
Let $\psi:\Gamma\to {\bf H}_{\Gamma}$ be the $k$-algebraic hull of $\Gamma$ and $\bf G$ the Zariski-closure of $\psi(\Delta)$ in  ${\bf H}_{\Gamma}$.
Then the algebraic group $\bf G$ is the $k$-algebraic hull of $\Delta$ and we have  ${\bf U}_{\Delta}={\bf U}_{\Gamma}$.
\end{lemma}
\begin{proof}

Let ${\bf H}_{\Gamma}^{0}$ be the identity component of ${\bf H}_{\Gamma}$.
Since $\bf G$ is a closed finite index subgroup of ${\bf H}_{\Gamma}$, we have ${\bf H}_{\Gamma}^{0}\subset {\bf G}$.
Since $\Gamma$ is virtually polycyclic, ${\bf H}_{\Gamma}^{0}$ is solvable.
Hence we have ${\bf U}({\bf H}_{\Gamma})=({\bf H}_{\Gamma}^{0})_{unip}={\bf U}({\bf G})$. 
Since ${\rm rank}\, \Gamma={\rm rank}\,\Delta$, we have
\[{\rm dim}\,{\bf U}({\bf G})={\rm rank}\,\Delta,
\]
and we have
\[ Z_{{\bf G}^{\prime}}({\bf U}({\bf G}))\subset Z_{{\bf H}_{\Gamma}}({\bf U}({\bf H}_{\Gamma}))\subset {\bf U}({\bf H}_{\Gamma})= {\bf  U}({\bf G}).
\]
Hence the lemma follows.
\end{proof}

\begin{lemma}\label{a-l-2}{\rm (\cite[Proof of Theorem 4.34]{R})}
Let $G$ be a simply connected solvable Lie group with a lattice $\Gamma$.
Let $\psi:G\to {\bf H}_{G}$ be the $\R$-algebraic hull of G and $\bf H^{\prime}$  the Zariski-closure of $\psi(\Gamma)$ in $\bf{H_{G}}$.
Then $\bf H^{\prime}$ is the $\R$-algebraic hull of $\Gamma$ and we have
${\bf U}_{G}={\bf U}_{\Gamma}$.
\end{lemma}

\section{Cohomology computations of aspherical manifolds with virtually torsion-free polycyclic fundamental groups  }\label{cohoo}
 Let $\Gamma$ be a torsion-free virtually polycyclic group and $\bf{H_{\Gamma}}$  the $\Q$-algebraic hull of $\Gamma$.
Denote $H_{\Gamma}=\bf{H_{\Gamma}}(\R)$. 
Let $U_{\Gamma}$ be the unipotent radical of $H_{\Gamma}$ and let $T$ be a maximal reductive subgroup.
Then $H_{\Gamma}$ decomposes as a semi-direct product $H_{\Gamma}=T\ltimes U_{\Gamma} $.
Let $\frak{u}$ be the Lie algebra of $U_{\Gamma}$. Since the exponential map ${\exp}:{\frak u} \longrightarrow U_{\Gamma}$ is a diffeomorphism, $U_{\Gamma}$ is  diffeomorphic to $\R^n$ such that $n={\rm rank}\,\Gamma$.
The splitting $H_{\Gamma}=T\ltimes U_{\Gamma}$ gives rise to the affine action $\alpha :H_{\Gamma}\longrightarrow {\rm Aut}(U_{\Gamma})\ltimes U_{\Gamma} $ such that $\alpha $ is an injective homomorphism.

In \cite{B} Baues constructed a compact  aspherical manifold $M_{\Gamma}=\alpha(\Gamma)\backslash U_{\Gamma}$ with $\pi_{1}(M_{\Gamma})=\Gamma$.
We call $M_{\Gamma}$ a standard $\Gamma$-manifold.
\begin{theorem}{\rm (\cite[Theorem 1.2]{B})}\label{a-t-3}
Standard $\Gamma$-manifold is unique up to diffeomorphism.
\end{theorem}
Let $A^{\ast}(M_{\Gamma})$ be the de Rham complex of $M_{\Gamma}$.
Then $A^{\ast}(M_{\Gamma}) $ is  the set of   the $\Gamma$-invariant differential forms ${A^{\ast}(U_{\Gamma})}^{\Gamma}$ on $U_{\Gamma}$. 
Let $(\bigwedge \frak{u} ^{\ast})^{T}$ be the left-invariant forms on $U_{\Gamma}$ which are fixed by $T$.
Since $\Gamma\subset H_{\Gamma}=T\ltimes U_{\Gamma}$, we have the inclusion
\[(\bigwedge {\frak u} ^{\ast})^{T} ={A^{\ast}(U_{\Gamma})}^{H_{\Gamma}} \subset {A^{\ast}(U_{\Gamma})}^{\Gamma}= A^{\ast}(M_{\Gamma}).\]
\begin{theorem}{\rm (\cite[Theorem 1.8]{B})}\label{a-t-4}
This inclusion induces a cohomology isomorphism.
\end{theorem}

\section{Proof of Theorem \ref{thoo}}
\subsection{\bf The embedings of   solvable Lie algebras in  splittable Lie algebras}
The idea of this subsection is based on \cite{Re}.
Let $\g$ be a solvable Lie algebra and $\n=\{X\in \g\vert {\rm ad}_{X}\, \, {\rm is\,\, nilpotent}\}$.
Then $\n$ is the maximal nilpotent ideal of $\g$ and called the nilradical of $\g$.
 \begin{lemma}\label{nilnil}{\rm (\cite[p.58]{OV2})}
 We have $[\g, \g]\subset \n$.
 \end{lemma}
 Let $D(\g)$ denote the space of the derivations of $\g$.
 By the Jordan decomposition, we have the decomposition ${\rm ad}_{X}=d_{X}+n_{X}$ such that $d_{X} $ is a semi-simple operator and $n_{X}$ is  a nilpotent operator.
 \begin{lemma}{\rm (\cite[Proposition 3]{Re})}
 We have $d_{X}$, $n_{X}$ $\in D(\g)$.
 \end{lemma}
 Then we have the homomorphism $f:\g\to D(\g)$ such that $f(X)=d_{X}$ for $X\in\g$.
 Since  ${\rm ker}f=\n$, we have ${\rm Im} f\cong \g/\n$.
 
 Let $\bar{\g} ={\rm Im} f\ltimes\g$ and
 $\bar{\n}=\{X-d_{X}\in  \bar{\g}  \vert X\in\g\}$.
Since ${\rm ad}_{X-d_{X}}={\rm ad}_{X}-d_{X}$ on $\g$,  ${\rm ad}_{X-d_{X}}$ is a nilpotent operator.
So $\bar{\n}$ consists of nilpotent elements. 
\begin{proposition}\label{spli}
We have $d_{X}(\bar{\n})\subset \n$ for any $X\in\g$, $\bar{\n}$ is a nilpotent ideal of $\bar{\g}$ and $\bar{\g}={\rm Im} f\ltimes\bar{\n}$.
\end{proposition}
\begin{proof}
 
 By  Lie's theorem, we have a basis $X_{1},\dots,X_{l}$ of $\n\otimes \C$ such that ${\rm ad}_{\g}$ on $\n$ are represented  by upper triangular matrices.
 Then for any $X\in\g$, we have \\
 ${\rm ad}_{X}(X_{1})=a_{X,1}X_{1},  $\\
 ${\rm ad}_{X}(X_{2})=a_{X,2}X_{2}+b_{X,12}X_{1},$\\
 $\ \ \ \ \vdots$\\
 ${\rm ad}_{X}(X_{l})=a_{X,l}X_{l}+b_{X,l-1l}X_{l-1}+\dots +b_{X,1l}X_{1}.$\\
 We take $X_{l+1},\dots ,X_{l+m}$ such that  $X_{1},\dots, X_{l}, X_{l+1},\dots , X_{l+m}$ is a   basis of $\g\otimes\C$.
 By Lemma \ref{nilnil},  we have ${\rm ad}_{X}(X_{i})\in \n$.
 Hence we have \\
 ${\rm ad}_{X}(X_{l+1})=b_{X,ll+1}X_{l}+\dots +b_{X,1l+1}X_{1}$\\
 $\ \ \ \  \vdots$\\
  ${\rm ad}_{X}(X_{l+m})=b_{X,ll+m}X_{l}+\dots +b_{X,1l+m}X_{1}.$\\
  Then we have 
  \[ d_{X}(X_{i})=a_{X,i}X_{i}, \ \ \ \ \  \ \ \ 1\le i \le l,\]
 \[ \ \  \ \ \ \     \ d_{X}(X_{i})=0, \ \ \ \ \ \ \ \  \ \ \ l+1\le i \le l+m.\]
Hence we have $d_{X}(\g)\subset \n$ and $d_{X}(\bar{\n})\subset\n$.
This implies $[\bar{\g}, \bar{\g}]\subset\n$.
In particular, $\bar{\n}$ is an ideal of $\bar{\g}$.
Since $\bar{\n}$ consists of nilpotent elements,  $\bar{\n}$ is a nilpotent ideal.
By $\bar{\g}=\{d_{X}+Y-d_{Y}\vert X, Y\in \g\}$, we have  $\bar{\g}={\rm Im}f\ltimes\bar{\n}$.
\end{proof}
By this proposition,  we have the inclusion $i:\g\to D(\bar{\n})\ltimes \bar{\n}$ given by $i(X)=d_{X}+X-d_{X}$ for $X\in\g$.

\subsection{\bf Constructions of algebraic hulls of simply connected solvable Lie groups }
Let $G$ be a simply connected solvable Lie group and $\g$  the Lie algebra of $G$.
Let $N$ be the maximal normal nilpotent subgroup of $G$ which corresponds to the nilradical $\n$ of $\g$.
Consider the injection $i\colon \g\to  {\rm Im}f\ltimes\bar{\n}\subset D(\bar{\n})\ltimes \bar{\n}$ constructed in  the last subsection.
Let $\bar{N}$ be the simply connected Lie group which corresponds to $\bar{\n}$. Since the Lie algebra of ${\rm Aut}(\bar{ N})\ltimes\bar{ N}$ is $D(\bar{\n})\ltimes \bar{\n}$,
we have the Lie group homomorphism $I:G\to {\rm Aut}(\bar{N})\ltimes\bar{N}$  induced by the injective homomorphism $i:\g\to D(\bar{\n})\ltimes \bar{\n}$.
\begin{lemma}\label{injj}
The homomorphism $I:G\to {\rm Aut}(\bar{N})\ltimes\bar{N}$ is injective.
\end{lemma}
\begin{proof}
Since the restriction of  $i:\g\to D(\bar{\n})\ltimes \bar{\n}$ on $\n$ is injective, the restriction $I:G\to {\rm Aut}(\bar{N})\ltimes\bar{N}$ on $N$ is also injective.
Let $T_{f}$ be the subgroup of  ${\rm Aut} (\bar{N})$ which corresponds to ${\rm Im} f$.
We have $I:G\to T_{f}\ltimes\bar{N}$.
By    Proposition \ref{spli}, $\bar{\g}/\n={\rm Im} f\oplus \bar{\n}/\n$.
So we have the induced map $I:G/N\to T_{f}\times\bar{N}/N$ and it is sufficient to show that this map is injective.
Let $j:{\rm Im}f\oplus \bar{\n}/\n\to\bar{\n}/\n$ be the projection and $J:T_{f}\times\bar{N}/N\to \bar{N}/N$ be the homomorphism which corresponds to  $j$.
Since the composition
\[ j\circ i (X \ \ \ {\rm mod}\,\n)=X-d_{X} \ \ \  {\rm mod}\,\n\]  is surjective,  $j\circ i:\g/\n\to\bar{\n}/\n$ is an isomorphism.
Since $G/N$ and $\bar{N}/N$ are simply connected abelian groups, $J\circ I :G/N \to \bar{N}/N$ is also an isomorphism.
Hence $I:G/N\to T_{f}\times\bar{N}/N$ is injective.
 
\end{proof}

A simply connected nilpotent Lie group is considered as the real points of a unipotent $\R$-algebraic group(see \cite[p. 43]{OV}) by the exponential map.
We have the unipotent $\R$-algebraic group $\bf \bar{N}$ with ${\bf \bar{N}}(\R)=\bar{N}$.
We identify the group $\Aut_{a}({\bf \bar{N}})$ of automorphisms of algebraic groups with $\Aut(\n_{\C})$ and $\Aut_{a}({\bf \bar{N}})$ has the $\R$-algebraic group structure with ${\rm Aut}_{a}({\bf \bar{N}})(\R)= {\rm Aut} (N)$.
So we have the $\R$-algebraic group $ \Aut_{a} (\bf\bar{N})\ltimes \bar{N}$.
By the above lemma, we have the injection $I:G\to {\rm Aut} (N)\ltimes N =\Aut_{a} (\bf \bar{N})\ltimes \bar{N}(\R)$.
Let $\bf G$ be the Zariski-closure of $I(G)$ in $ \Aut_{a} (\bf  \bar{N})\ltimes \bar{N}$.

\begin{lemma}\label{u=n}
We have $\U(\bf G)=\bar{N}$.
\end{lemma}
\begin{proof}
Let $\bf T$ be the Zariski-closure of $T_{f}$ in ${\rm Aut}_{a}( \bar{\bf N})$.
Then ${\bf G }\subset {\bf T} \ltimes \bar{\bf N}$.
Since $\bf G$ is connected solvable and $\bf T$ consists of semi-simple automorphisms, we have ${\bf U}({\bf G})={\bf G}\cap  \bar{\bf N}$.
By this, it is sufficient to show $\dim {\bf U}({\bf G})=\dim \bar{\bf N}$.
Let $\bf N$ be the Zariski-closure of $I(N)$.
By $I(N)\subset \bar{N}$, we have  $ {\bf U}({\bf G})/{\bf N}={\bf U}({\bf G}/{\bf N})$.
Thus it is sufficient to show $ {\bf U} ({\bf G}/{\bf N})= G/N$.
Consider the induced map $I:G/N\to T_{f} \times \bar{N}/N$ as the proof of Lemma \ref{injj}.
The Zariski-closure of $I(G/N)$ in ${\bf T} \times \bar{\bf N}/{\bf N}$ is ${\bf G}/{\bf N}$.
Since ${\bf T} \times \bar{\bf N}/{\bf N}$ is commutative, the projection ${\bf T} \times \bar{\bf N}/{\bf N}\to  \bar{\bf N}/{\bf N}$ is an $\R$-algebraic group homomorphism.
Since we showed that $J\circ I:G/N \to \bar{N}/N$ is isomorphism In the proof of Lemma \ref{injj}, 
the image $J\circ I(G/N)$ is Zariski-dense in $\bar{\bf N}/{\bf N}$.
This implies $\bar{\bf N}/{\bf N}={\bf U}({\bf G}/{\bf N})$.
Hence the lemma follows.
\end{proof}
By this lemma we have the following proposition.
\begin{proposition}\label{hulll}
$\bf G$ is the algebraic hull of $G$ and the Lie algebra of the unipotent hull ${\bf U}_{G}$ is $\bar{\n}_{\C}$.
\end{proposition}
\begin{proof}
We show that $\bf G$ satisfies the properties of the algebraic hull of $G$.
We have $\dim\U({\bf G})=\dim\bar{ {\bf N}}=\dim G$.
Let $(t,\ x)\in Z_{\bf G}(\U({\bf G}))\subset \Aut_{a} \bf \bar{N}\ltimes \bar{N}$.
Since $\U(\bf G)=\bf N$ and t is a semi-simple automorphism, we have $t(y)=y$ for any $y\in \bar{\bf N}$.
So we have $t={\rm id}_{\bf \bar{N}}$.
We have  $Z_{\bf G}(\U({\bf G}))\subset \U({\bf G})$.
Hence the proposition follows.
 
\end{proof}
\subsection{\bf Proof of Theorem \ref{thoo}}
We first prove:
\begin{theorem}\label{abab}
Let $G$ be a simply connected solvable Lie group.
Then ${\bf U}_{G}$ is abelian if and only if  $G=\R^{n}\ltimes_{\phi} \R^{m}$ such that the action $\phi:\R^{n}\to {\rm  Aut} (\R^{m})$ is semi-simple.
\end{theorem}
\begin{proof}
Consider the inclusion $i:\g \to {\rm Im}f\ltimes \bar{\n}$.
By the above argument, the Lie algebra of ${\bf  U}_{G}$ is $\bar{\n}_{\C}$.
Suppose $G=\R^{n}\ltimes_{\phi} \R^{m}$ such that the action $\phi:\R^{n}\to {\rm  Aut}\, \R^{m}$ is semi-simple.
It is sufficient to show  $\bar{\n}=\{X-d_{X}\vert X\in \g\}\subset {\rm Im}f\ltimes \g$ is an abelian Lie algebra. 
Let $X, Y\in \g$ and $X=X_{1}+X_{2} $, $Y=Y_{1}+Y_{2}$ be the decompositions induced by the semi-direct product $\g=\R^{n}\ltimes_{\phi_{\ast}} \R^{m}$.
Then we have $d_{X_{2}}=0$, $d_{Y_{2}}=0$, $[X_{1}, Y_{1}]=0$ and $[X_{2}, Y_{2}]=0$
by the assumption.
Hence we have 
\[[X-d_{X}, Y-d_{Y}]=[X_{1},Y_{2}]+[X_{2},Y_{1}]-d_{X_{1}}(Y_{2})+d_{Y_{1}}(X_{2}).
\]
Since the action $\phi_{\ast} $ is semi-simple, we have $d_{X_{1}}(Y_{2})=[X_{1},Y_{2}]$ and  $d_{Y_{1}}(X_{2})=[Y_{1},X_{2}]$.
Therefore we have $[X-d_{X}, Y-d_{Y}]=0$. This implies $\bar{\n}$ is abelian.

Conversely we assume ${\bf U}_{G}$ is abelian.
By Proposition \ref{hulll}, $\bar{\n}$ is abelian.
By $[\g,\g]\subset \n$, $\g$ is two-step solvable.
By \cite[Lemma 4.1]{Bur}, we have the decomposition $\g= {\frak a}\ltimes \g^{\infty}$
for some nilpotent subalgebra $\frak a$ of $\g$ where $\g^{\infty}=\bigcap \g^{i}$ for the lower central series $\g=\g^{0}\supset \g^{1}\supset  \g^{2}\supset\dots$ of $\g$.
Since $\bar{\n}$ is abelian, the  subspace $\{X-d_{X}\vert X\in{\frak a}\}$ is a abelian subalgebra of $\bar \n$.
Since $\frak a$ is nilpotent,  the  Lie algebra $\{X-d_{X}\vert X\in{\frak a}\}$ is  identified with $\frak a$.
Hence $\frak a$ is abelian.
Finally we show that the action of $\frak a$ on $\g^{\infty}$ is semi-simple.
We suppose that ${\rm ad}_{X}$ on $\g^{\infty}$ is not semi-simple for some $X\in \frak a$.
Then the action of ${\rm ad}_{X}-d_{X}$ on $\g^{\infty}$ is non-trivial.
Since we have $\bar{\n}=\{X-d_{X}\vert X\in \g\}\subset {\rm Im}f\ltimes \bar{\n}$, we have $[\bar{\n},{\frak a}]\not= \{0\}$.
This contradicts $\bar{\n}$ is abelian. 
Hence the action of $\frak a$ on $\g^{\infty}$ is semi-simple and we have the theorem.
\end{proof}

\begin{proof}[{\bf Proof of Theorem \ref{thoo}}]
By Theorem \ref{subpo}, we have a finite index subgroup of $\Gamma$ which is isomorphic to a lattice of some simply connected solvable Lie  group $G$.
By Lemma \ref{a-l-1} and \ref{a-l-2}, we have ${\bf U}_{\Gamma}={\bf U}_{G}$.
Hence by Theorem \ref{abab} we have the theorem.
\end{proof}

\begin{remark}
A virtually polycyclic group $\Gamma$ has the maximal nilpotent normal subgroup called the nilradical of $\Gamma$.
Since the nilradical of $\Gamma$ is contained in ${\bf U}_{\Gamma}$ (see \cite[Proposition A.7]{B}), if ${\bf U}_{\Gamma}$ is abelian then the nilradical of $\Gamma$ is also abelian.
But the converse is not true.
Consider $G= \R\ltimes _{\phi}\R^{4}$ with
\[\phi(t)= \left(
\begin{array}{cccc}
e^{rt}&0&0&0\\
0&e^{-rt}&0&0\\
0&0&1&t\\
0&0&0&1
\end{array}
\right).
\]
Then for some $r\not=0$ $G$ has a lattice $\Z\ltimes_{\phi} \Gamma^{\prime\prime}$ for a lattice $\Gamma^{\prime\prime}$ of $\R^{4}$.
We have ${\bf U}_{\Gamma}={\bf U}_{G}=\C^{2}\times {\bf U}_{3}(\C)$ and it is not abelian.
On the other hand the nilradical of $\Gamma$(resp. $G$) is isomorphic to $\Z^{4}$(resp. $\R^{4})$. 
\end{remark}

\section{Formality and  hard Lefschetz properties  of aspherical manifolds}  
\subsection{\bf Formality} 
We review the definition of formality and  prove Proposition \ref{abu}.
\begin{definition} A {\em differential graded algebra} (called  {\em DGA}) is a graded $\R $-algebra $A^{\ast}$  with the following properties: \\
(1)
$A^{\ast}$ is graded commutative, i.e.
\[y\wedge x=(-1)^{p\cdot q}x\wedge y \ \ \ x\in A^{p} \ \ \  y\in A^{q}.
\]
(2)
There is a differential operator $d:A\rightarrow A$ of degree one such that $d\circ d=0$ and
\[d(x\wedge y)=dx\wedge y+(-1)^{p}x\wedge dy  \ \ \  x\in A^{p}.
\]
\end{definition}
Let $A$ and $B$ be DGAs.
If  a morphism of graded algebra $\varphi : A \rightarrow B$ satisfies $d\circ \varphi =\varphi \circ d$, we call $\varphi $ a morphism of DGAs.
If a morphism of DGAs induces a cohomology isomorphism, we call it a quasi-isomorphism.
\begin{definition}
A and B are {\em weakly equivalent} if there is a finite diagram of DGAs 
\[A\leftarrow C_{1}\rightarrow C_{2}\leftarrow\cdot \cdot \cdot \leftarrow C_{n}\rightarrow B
\]
such that all the morphisms are quasi-isomorphisms.
\end{definition}
Let $M$ be a smooth manifold.
The de Rham complex  $A^{\ast}(M)$ of $M$ is  a DGA.
The cohomology algebra $H^{\ast}(M,\R)$ is a DGA with $d=0$.
\begin{definition}\label{df}
A smooth manifold $M$ is {\em formal} if $A^{\ast}(M)$ and $H^{\ast}(M,\R)$ are weakly equivalent.
\end{definition}
\begin{proposition}\label{au}
Let $\Gamma$ be a  torsion-free virtually polycyclic group.
If the  unipotent hull  $U_{\Gamma}$ is abelian,  the standard $\Gamma$-manifold $M_{\Gamma}$ is formal.
\end{proposition}
\begin{proof}
We use same notations as in Section \ref{cohoo}.
If the  $k$-unipotent hull of $\Gamma$ is abelian, $(\bigwedge {\frak u} ^{\ast},d)=(\bigwedge {\frak u} ^{\ast},0)$.
By  Theorem \ref{a-t-4}, we have the diagram of  DGAs 

\[A^{\ast}(M_{\Gamma})\leftarrow ((\bigwedge {\frak u} ^{\ast})^{T})=H^{\ast}(M_{\Gamma})
\]
such that the map $A^{\ast}(M_{\Gamma})\leftarrow ((\bigwedge {\frak u} ^{\ast})^{T})$  is a quasi-isomorphism.
Hence the  proposition follows.
 
\end{proof}
Hence we have Proposition \ref{abu}.

\subsection{\bf The hard Lefschetz property}
We review the definition of the hard Lefschetz property and  prove Proposition \ref{LEE}.
\begin{definition}\label{dhd}
Let $(M, \omega)$ be a $2n$-dimensional symplectic manifold.
We say that $(M,\omega)$ satisfies the {\em hard Lefschetz property} if the linear map
\[[\omega^{n-i}]\wedge : H^{i}(M, \R )\to H^{2n-i}(M, \R)\]
is an isomorphism for any $0\leq i\leq n$.
\end{definition}
\begin{proof}[{\bf Proof of Proposition \ref{LEE}}]
As in the proof of Proposition \ref{abu}, we have an isomorphism $(\bigwedge {\frak u}^{\ast})^{T}\cong H^{\ast}(M, \R)$.
Consider the cohomology class of a symplectic form $\omega$ on $M$.
 We have $\omega_{0}\in (\bigwedge^{2} {\frak u}^{\ast})^{T}$ which represents the cohomology class $[\omega]\in H^{2}(M,\R)$.
Since $\omega^{n}_{0}\not=0$ for $2n=\dim \frak u=\dim M$, $\omega_{0}$ is a symplectic form on the vector space $\frak u$.
Since the linear map
\[\omega_{0}^{n-i}\wedge : \bigwedge {\frak u}^{i}\to \bigwedge {\frak u}^{2n-i}\]
is injective for any $0 \leq i\leq n$ by the hard Lefschetz property of a torus, the restriction
\[\omega_{0}^{n-i}\wedge : (\bigwedge {\frak u}^{i})^{T}\to (\bigwedge {\frak u}^{2n-i})^{T}\]
is also injective and so
\[[\omega^{n-i}]\wedge : H^{i}(M, \R)\to H^{2n-i}(M, \R)\]
is injective and thus it is  an isomorphism by the Poincar\'e duality.
Hence we have the proposition.
\end{proof}

\section{Examples}

\begin{example}\label{1}
Let $G=\R\ltimes_{\phi} \R^{2}$ with $\phi(t)= \left(
\begin{array}{ccc}
e^{rt}&0\\
0&e^{-rt}
\end{array}
\right)$.
Then for some $r\not=0$ $\phi(1)$ is conjugate to an element of $SL_{2}(\Z)$.
Hence we have a lattice $\Gamma=\Z\ltimes \Z^{2}$.
$G\times \R$ has a left-invariant symplectic form.
In \cite{FG}(see also \cite{T}) by direct computations  Fernandez and Gray showed that $G/\Gamma\times S^{1}$ is formal and satisfie the hard Lefschetz property and admits no Complex structure.
This is also a simple example for the result of this paper.
\end{example}

\begin{example}\label{cp}

Let $G=\C\ltimes_{\phi} \C^{2}$ with $\phi(x)= \left(
\begin{array}{ccc}
e^{x}&0\\
0&e^{-x}
\end{array}
\right)$.
Then the cochain complex $(\bigwedge \g^{\ast}, d)$ of the Lie algebra of $G$ 
is given by:\\
\[\g^{\ast}=\langle x_{1}, x_{2}, y_{1}, y_{2}, z_{1}, z_{2} \rangle ,
\]
\[dx_{1}=dx_{2}=0 ,
\]
\[dy_{1}=-x_{1}\wedge y_{1}+x_{2}\wedge y_{2},\, dy_{2}=-x_{2}\wedge y_{1}
-x_{1}\wedge y_{2}, 
\]
\[dz_{1}=x_{1}\wedge z_{1}-x_{2}\wedge z_{2},\, dz_{2}=x_{1}\wedge z_{2}+x
_{2}\wedge z_{1}.
\]
We have an invariant symplectic form $\omega=x_{1}\wedge x_{2}+z_{1}\wedge y_{1} +y_{2}\wedge z_{2} $.
For some $p,q\in \R$ $\phi(p\Z+\sqrt{-1}q\Z)$ is conjugate to a subgroup of $SL_{4}(\Z)$ and hence we have a lattice $\Gamma=(p\Z+\sqrt{-1}q\Z)\ltimes \Gamma^{\prime\prime}$ for  a lattice $\Gamma^{\prime\prime}$ of $\C^{2}$(see \cite{Na} and \cite{Hd}).
For any lattice $\Gamma$, $G/\Gamma$ is complex, symplectic with the  hard Lefschetz property and formal but not K\"ahler.
 
\begin{remark}
For a Lie group $G$ in Example \ref{cp}, the de Rham cohomology of $G/\Gamma$ depends on a choise  
of a lattice $\Gamma$.
Under some conditions, the de Rham cohomology of a solvmanifold $G/\Gamma $ is isomorphic to the cohomology of Lie algebra $\g$ of $G$ (see \cite{Hatt}, \cite[Section 7]{R}).
But for a general solvmanifold $G/\Gamma$ it is difficult to compute the de Rham cohomology of $G/\Gamma$.
By the results of this paper, for a Lie group  $G=\R^{n}\ltimes _{\phi}\R^{m}$ with the semi-simple action $\phi$, we can say that $G/\Gamma$ is formal and hard  Lefschetz for any lattice $\Gamma$ even if an isomorphism $H^{\ast}(G/\Gamma,\R)\cong H^{\ast}(\g)$ fails to hold.
\end{remark}
\end{example}
\begin{example}
Let $G=\R\ltimes_{\phi}\R^{4}$ with \[\phi(t)=\left(
\begin{array}{cccc}
e^{pt}\cos(qt)&-e^{pt}\sin(qt)&0&0\\
e^{pt}\sin(qt)&e^{pt}\cos(qt)&0&0\\
0&0&e^{-pt}\cos(-qt) &-e^{-pt}\sin(-qt) \\
0&0&e^{-pt}\sin(-qt) &e^{-pt}\cos(-qt) 
\end{array}
\right).
\]
Then for $p, q$ as Example \ref{cp}, $\phi(1)$ is conjugate to an element of $SL_{4}(\Z)$ and hence $G$ has a lattice $\Gamma=\Z\ltimes \Gamma^{\prime\prime}$ for  a lattice $\Gamma^{\prime\prime}$ of $\R^{4}$.
The cochain complex $(\bigwedge (\g\oplus\R)^{\ast}, d)$ of the Lie algebra of $G\times \R$ 
is given by:
\[(\g\oplus\R)^{\ast}=\langle w,x_{1},x_{2},x_{3},x_{4}, y\rangle,\]
\[dx_{1}=-pw\wedge x_{1}+qw\wedge x_{2},\, \,  dx_{2}=-qw\wedge x_{1}-pw\wedge x_{2},\]
\[dx_{3}=pw\wedge x_{3}-qw\wedge x_{4},\,\, dx_{4}=qw\wedge x_{3}+pw\wedge x_{4}.\]
We have a left-invariant symplectic form $\omega=w\wedge y+x_{1}\wedge x_{3}+x_{4}\wedge x_{2}$.
We regard $w+\sqrt{-1}y, x_{1}+\sqrt{-1}x_{2}, x_{3}+\sqrt{-1} x_{4}$ as $(1,0)$-forms, we obtain a left-invariant complex structure.
By the result of this paper, for any lattice $\Gamma$, $G/\Gamma\times S^{1}$ is formal and  any symplectic form on  $G/\Gamma\times S^{1}$ satisfies the hard Lefschetz property.
\begin{remark}
In \cite{Boc}, Bock studies formality and the hard Lefschetz property of solvmanifolds of dimension$\le 6$ by direct computations.
The cohomology of $G/\Gamma$ may vary for a choice of $\Gamma$ and  Bock does not decide whether $G/\Gamma\times S^{1}$ is formal and satisfies the hard Lefschetz property.

\end{remark}
\end{example}
By combining the above examples we obtain:
\begin{example}
Let $G=\R^{2}\ltimes_{\phi} \R^{2k+4(l+m+n)}$ such that 
\begin{multline*}
\phi(t_{1},t_{2})=\bigoplus _{i=1}^{k}\left(
\begin{array}{ccc}
\cos a_{i}t_{1}&-\sin a_{i}t_{1}\\
\sin a_{i}t_{1}&\cos a_{i}t_{1}
\end{array}
\right)\oplus\\
\bigoplus_{i=1}^{l} \left(
\begin{array}{cccc}
e^{b_{i}t_{1}}&0&0&0\\
0&e^{-b_{i}t_{1}}&0&0\\
0&0&e^{b_{i}t_{1}}&0\\
0&0&0&e^{-b_{i}t_{1}}
\end{array}
\right)\oplus  \\
\bigoplus_{i=1}^{m} \left(
\begin{array}{cccc}
e^{c_{i}t_{1}}\cos(d_{i}t_{2})&-e^{c_{i}t_{1}}\sin(d_{i}t_{2})&0&0\\
e^{c_{i}t_{1}}\sin(d_{i}t_{2})&e^{c_{i}t_{1}}\cos(d_{i}t_{2})&0&0\\
0&0&e^{-c_{i}t_{1}}\cos(-d_{i}t_{2})&-e^{-c_{i}t_{1}}\sin(-d_{i}t_{2}) \\
0&0&e^{-c_{i}t_{1}}\sin(-d_{i}t_{2})&e^{-c_{i}t_{1}}\cos(-d_{i}t_{2}) 
\end{array}
\right)\oplus\\
\bigoplus_{i=1}^{n} \left(
\begin{array}{cccc}
e^{e_{i}t_{1}}\cos(f_{i}t_{1})&-e^{e_{i}t_{1}}\sin(f_{i}t_{1})&0&0\\
e^{e_{i}t_{1}}\sin(f_{i}t_{1})&e^{e_{i}t_{1}}\cos(f_{i}t_{1})&0&0\\
0&0&e^{-e_{i}t_{1}}\cos(-f_{i}t_{1})&-e^{-e_{i}t_{1}}\sin(-f_{i}t_{1}) \\
0&0&e^{-e_{i}t_{1}}\sin(-f_{i}t_{1})&e^{-e_{i}t_{1}}\cos(-f_{i}t_{1}) 
\end{array}
\right).
\end{multline*}
We write $A\oplus B=\left(
\begin{array}{ccc}
A&0\\
0&B
\end{array}
\right)$ for matrices $A,B$.

We suppose $a_{i}=\frac{2\pi}{K_{i}}$ for $K_{i}=2,3,4,$ or $6$, $b_{i}=rL_{i}$ for $r$ as Example \ref{1} and $L_{i}\in\Z$, $c_{i}=pM_{i}$, $d_{i}=qM_{i}^{\prime}$, $e_{i}=pN_{i}$ and $f_{i}=qN^{\prime}_{i}$ for $p,q$ as Example \ref{cp} and $M_{i},M_{i}^{\prime}, N_{i}, N_{i}^{\prime}\in \Z$.
Then  each component of $\phi(\Z^{2})$ for the direct product is conjugate to a subgroup of   $SL_{2}(\Z)$ or $SL_{4}(\Z)$ and hence we have a lattice $\Gamma=\Z^{2}\ltimes \Gamma^{\prime\prime}$ for a lattice $\Gamma^{\prime\prime}$ of $ \R^{2k+4(l+m+n)}$.
 the cochain complex $(\bigwedge \g^{\ast}, d)$ of the Lie algebra of $G$ 
is given by:
\[\g^{\ast}=\langle u_{1} ,u_{2}, w_{1},\dots,w_{2k}, x_{1},\dots ,x_{4l}, y_{1},\dots ,y_{4m}, z_{1},\dots ,z_{4n}  \rangle ,
\]
\[ d u_{1} =d u_{2} =0,\]
\[dw_{2i-1}=a_{i}u_{1} \wedge w_{2i},\, \, dw_{2i}=-a_{i}u_{1} \wedge w_{2i-1}, \ \ \ \ (1\leq i\leq  k),
\]
\[dx_{2i-1}=-b_{i}u_{1} \wedge x_{2i-1},\, \, dx_{2i}= b_{i}u_{1}\wedge x_{2i}, \ \ \ \ (1\leq i\leq 2l),
\]
\[dy_{4i-3}=-c_{i}u_{1}\wedge y_{4i-3}+d_{i}u_{2}\wedge y_{4i-2},\, \, dy_{4i-2}=-d_{i}u_{2}\wedge y_{4i-3}-c_{i}u_{1}\wedge y_{4i-2},
\]
\[dy_{4i-1}=c_{i}u_{1}\wedge y_{4i-1}-d_{i}u_{2}\wedge y_{4i},\, \, dy_{4i}=d_{i}u_{2}\wedge y_{4i-1}+c_{i}u_{1}\wedge y_{4i},
\]
\[(1\le i\le m),\]
 \[dz_{4i-3}=-e_{i}u_{1}\wedge z_{4i-3}+f_{i}u_{1}\wedge z_{4i-2},\, \, dz_{4i-2}=-f_{i}u_{1}\wedge z_{4i-3}-e_{i}u_{1}\wedge z_{4i-2},
\]
\[dz_{4i-1}=e_{i}u_{1}\wedge z_{4i-1}-f_{i}u_{1}\wedge z_{4i},\, \, dz_{4i}=f_{i}u_{1}\wedge z_{4i-1}+e_{i}u_{1}\wedge z_{4i},\]
\[(1\le i\le n).\]
$G$ has a left-invariant symplectic form
\begin{multline*}
\omega 
=u_{1}\wedge u_{2}+\sum_{i=1}^{k}w_{2i-1}\wedge w_{2i}+\sum _{i=1}^{2l}x_{2i-1}\wedge x_{2i}\\
+\sum_{i=1}^{m}(y_{4i-3}\wedge y_{4i-1}+y_{4i}\wedge y_{4i-2})+\sum_{i=1}^{n}(z_{4i-3}\wedge z_{4i-1}+z_{4i}\wedge z_{4i-2}).
\end{multline*}
Regarding 

\[u_{1}+\sqrt{-1}u_{2}, \]
\[w_{2i-1}+\sqrt{-1}w_{2i}\, \, (1\le i\le k),\]
\[ x_{4i-3}+\sqrt{-1}x_{4i-1}, x_{4i-2}+\sqrt{-1}x_{4i}\, \, (1\le i\le l),\] \[y_{2i-1}+\sqrt{-1}y_{2i}\,\, (1\le i\le 2m), \]
\[z_{2i-1}+\sqrt{-1}z_{2i}\, \, (1\le i\le 2n)
\]
 as $(1,0)$-forms, we have a left-invariant complex structure on $G$.
By the results of this paper, for any lattice $\Gamma$ of $G$, $G/\Gamma$ is formal and satisfies the hard Lefschetz property but admits no K\"ahler structure.
\end{example}

\begin{example}[Oeljeklaus-Toma manifolds]
We apply the result of this paper to non-K\"ahler complex manifolds constructed by Oeljeklaus and Toma in \cite{OT}.
Let $K$ be a finite extension field of $\Q$ with the degree $s+2t$ for positive integers $s,t$.
Suppose $K$ admits embeddings $\sigma_{1},\dots \sigma_{s},\sigma_{s+1},\dots, \sigma_{s+2t}$ into $\C$ such that $\sigma_{1},\dots ,\sigma_{s}$ are real embeddings and $\sigma_{s+1},\dots, \sigma_{s+2t}$ are complex ones satisfying $\sigma_{s+i}=\bar \sigma_{s+i+t}$ for $1\le i\le t$. 
We can choose $K$ admitting such embeddings(see \cite{OT}).
Denote ${\mathcal O}_{K}$ the ring of algebraic integers of $K$, ${\mathcal O}_{K}^{\ast}$ the group of units in ${\mathcal O}_{K}$ and 
\[{\mathcal O}_{K}^{\ast\, +}=\{a\in {\mathcal O}_{K}^{\ast}: \sigma_{i}>0 \,\, {\rm for \,\,  all}\,\, 1\le i\le s\}.
\]  
Define $l:{\mathcal O}_{K}^{\ast\, +}\to \R^{s+t}$ by 
\[l(a)=(\log \vert \sigma_{1}(a)\vert,\dots ,\log \vert \sigma_{s}(a)\vert , 2\log \vert \sigma_{s+1}(a)\vert,\dots ,2\log \vert \sigma_{s+t}(a)\vert)
\]
for $a\in {\mathcal O}_{K}^{\ast\, +}$.
Then by Dirichlet's units theorem, $l({\mathcal O}_{K}^{\ast\, +})$ is a lattice in the vector space $L=\{x\in \R^{s+t}\vert \sum_{i=1}^{s+t} x_{i}=0\}$.
For the projection $p:L\to \R^{s}$ given by the first $s$ coordinate functions.
Then we have a  subgroup $U$ with the rank $s$ of ${\mathcal O}_{K}^{\ast\, +}$ such that $p(l(U))$ is a lattice in $\R^{s}$.
We have the action of $U\ltimes{\mathcal O}_{K}$ on $H^{s}\times \C^{t}$
such that 
\begin{multline*}
(a,b)\cdot (x_{1}+\sqrt{-1}y_{1},\dots ,x_{s}+\sqrt{-1}y_{s}, z_{1},\dots ,z_{t})\\
=(\sigma_{1}(a)x_{1}+\sigma_{1}(b)+\sqrt{-1} \sigma_{1}(a)y_{1}, \dots ,\sigma_{s}(a)x_{s}+\sigma_{s}(b)+\sqrt{-1} \sigma_{s}(a)y_{s},\\
 \sigma_{s+1}(a)z_{1}+\sigma_{s+1}(b),\dots ,\sigma_{s+t}(a)z_{t}+\sigma_{s+t}(b)).
\end{multline*}
In \cite{OT} it is proved that the quotient $X(K,U)=H^{s}\times \C^{t}/U\ltimes{\mathcal O}_{K}$ is compact.
We call this complex manifold a  Oeljeklaus-Toma(OT) manifold with $(s,t)$.
By this construction we give  solvaminfold-presentations $G/\Gamma$ of  OT-manifolds with $(s,t)$.
We consider $p(l(U))\ltimes_{\phi}(\R^{s}\times \C^{t})$ with
\begin{multline*}
\phi(t_{1},\dots ,t_{s})\\
=\left(
\begin{array}{cccccc}
e^{t_{1}}&  &  && &\\
&     \ddots &   && & \\
&  &e^{t_{s}}&&&\\
&&& \sigma_{s+1}\circ \sigma_{1}^{-1}(e^{t_{1}})& &\\
 &&& &\ddots&\\
&&&&& \sigma_{s+t}\circ \sigma_{1}^{-1}(e^{t_{1}})
\end{array}
\right)
\end{multline*}
for $(t_{1},\dots ,t_{s})\in p(l(U))$.
Then for some lattice $\Gamma^{\prime\prime}$ of $\R^{s}\times \C^{t}$, we have  $p(l(U))\ltimes_{\phi}\Gamma^{\prime\prime}\cong U\ltimes{\mathcal O}_{K}$.
Since $p(l(U))$ is a lattice of $\R^{s}$, we have an extension of $\phi$ on $\R^{s}$ and $U\ltimes{\mathcal O}_{K}$ can be seen as a lattice of $\R^{s}\ltimes_{\phi}(\R^{s}\times \C^{t})$.
Thus OT-manifolds are formal complex solvmanifolds not admitting K\"ahler structure.
\begin{remark}
For $t=1$, OT-manifolds $X(K,U)$ admit LCK(locally conformal K\"ahler) structures.
\end{remark}
\begin{remark}
We call $X(K,U)$ simple type if the action of $U$ on $\mathcal O$ admits no proper non-trivial submodule of lower rank.
If  $X(K,U)$ is simple type, then in \cite{OT}  it is proved that the second Betti number is $b_{2}=\frac{s(s-1)}{2}$.
Then the second cohomology $H^{2}(X(K,U),\R)$ is spanned by $\{[dt_{i}\wedge dt_{j}]\}_{1\le i<j\le s}$ and hence simple type OT-manifolds admit no symplectic structure.
\end{remark}
\end{example}
\begin{example}
Infra-solvmanifolds appear in study of geometries of $3$-manifolds.
See \cite{PS} for the general theory of geometries of $3$-manifolds.
A compact aspherical  $3$-manifold $M$ with the virtually solvable fundamental group admits a  one of the three geometries $E^{3}, Nil, Sol$ i.e. $M$ is diffeomorphic to $G/\Gamma$ such that $G$ is $\R^{3}$, ${\bf U}_{3}(\R)$ or $\R\ltimes _{\phi}\R^{2}$ as an Example \ref{1} with a left-invariant metric  and $\Gamma\subset C\ltimes G$ is a lattice for the group $C$ of isometric automorphisms of $G$.
In the $E^{3}$ case, $\Gamma$ is virtually abelian  by Bieberbach's first theorem.
In the $Sol$ case, $C$ is finite(see \cite{PS}).
Hence a compact  $3$-manifold $M$  admitting the geometry $E^{3}$ or $Sol$ is formal.
\end{example}

\section{Remarks} 
In this Section we give an example of a formal standard $\Gamma$-manifold with the hard Lefschetz property such that $U_{\Gamma}$ is not 
abelian.
In addition this is also  an example of formal manifold satisfying the hard 
Lefschetz property such that it is finitely covered by a non-formal 
manifold  not satisfying the hard Lefschetz property.
We notice that compact manifolds finitely covered by non-K\"ahler 
manifolds are not K\"ahler.

Let $\Gamma=\Z \ltimes _{\phi}\Z^{2}$ such that 
 for $t\in\Z$
\[\phi(t)
=\left(
\begin{array}{cc}
(-1)^{t}& (-1)^{t} t  \\
0&     (-1)^{t}  
\end{array}
\right).
\]
\begin{lemma}
The algebraic hull of $\Gamma$ is given by $\bf H_{\Gamma}= \{\pm 1\} \ltimes {\bf U}_{3}(\C)$ such that 
 
\[ (-1)\cdot \left(
\begin{array}{ccc}
1&  x&  z  \\
0&     1&y      \\
0& 0&1 
\end{array}
\right)=
\left(
\begin{array}{ccc}
1&  x&  (-1)z  \\
0&     1&(-1)y      \\
0& 0&1 
\end{array}
\right).
\]

\end{lemma}
\begin{proof}
We have the inclusion
\[\Gamma \cong \left( (-1)^{x}, \left(
\begin{array}{ccc}
1&  x&  z  \\
0&     1&y      \\
0& 0&1 
\end{array}
\right)\right) \subset
 \{\pm 1\}  \ltimes{\bf U}_{3}(\C).
\]
Then $\Gamma$ is Zariski-dense in $\{\pm 1\}\ltimes {\bf U}_{3}(\C)$ and
${\rm rank}\, \Gamma=3=\dim{\bf U}_{3}(\C)$.
Since the action of $\{\pm 1\}$  on ${\bf U}_{3}(\C)$ is faithful,   the centralizer  of ${\bf U}_{3}(\C)$ is contained in  ${\bf U}_{3}(\C)$.
Hence the lemma follows.
\end{proof}
We have ${\bf H}_{\Gamma}(\R)=\{\pm 1\}\ltimes U_{\Gamma}$ such that  $U_{\Gamma}={\bf U}_{3}(\R)$.
Let $\mathfrak u$ be the Lie algebra of $U_{\Gamma} $.
We have $\mathfrak{u} =\langle X_{1}, X_{2}, X_{3} \rangle $ such that the bracket is given by 
\[ [X_{1},X_{2}]=-[X_{2},X_{1}]=X_{3}.
\]
The $\{\pm 1\}$-action on $\mathfrak u$ is given by 
\[ (-1)\cdot X_{1}=X_{1},\ (-1)\cdot X_{i}=-X_{i} \ \ \ \ \ \ \ \ i=2,3\]

Let  $x_{1}, x_{2}, x_{3} $ be the basis of $\mathfrak{u}^{\ast}$ which is dual to $X_{1},X_{2},X_{3}$.
Then the DGA $(\bigwedge {\mathfrak u} ^{\ast})^{\{\pm 1\}}$ 
is the subalgebra of $\bigwedge {\mathfrak u} ^{\ast}$ generated by $\{x_{1}, x_{2}\wedge x_{3}\}$ and the derivation on $(\bigwedge {\mathfrak u} ^{\ast})^{\{\pm 1\}}$ is trivial.
Let $M_{\Gamma}$ be the standard $\Gamma$-manifold.
Then by Theorem \ref{a-t-4}, we have the quasi-isomorphism $(\bigwedge {\mathfrak u} ^{\ast})^{\{\pm 1\}}\to A^{\ast}(M_{\Gamma})$. Since the derivation on $(\bigwedge {\mathfrak u} ^{\ast})^{\{\pm 1\}}$ is trivial, we have the isomorphism $(\bigwedge {\mathfrak u} ^{\ast})^{\{\pm 1\}}\cong H^{\ast}(M)$. Hence we have:
\begin{proposition}
$M_{\Gamma}$ is formal.
\end{proposition}

\begin{remark}
Since $U_{\Gamma}$ is not abelian, the converse of 
Proposition \ref{au} is not true. 
\end{remark}
\begin{remark}\label{fi} 
We have the finite index subgroup $2\Z\ltimes\Z^{2}$ which is nilpotent.
So $\Gamma$ is virtually  nilpotent but not virtually abelian.
By the result of \cite{H}, $K(2\Z\ltimes\Z^{2},1)$ is not formal. But for the finite extension group $\Gamma$, $K(\Gamma, 1)$ is formal.
\end{remark}
\begin{remark}
Since $\{\pm 1\}$ acts isometrically on  $U_{\Gamma}$ with the invariant metric, $M_{\Gamma}$ admits the $Nil $ geometry.
So we have a formal $3$-dimensional compact manifold admitting the $Nil$ geometry.
\end{remark} 

Let $\Delta=\Gamma\times \Z$. 
Then we have $H_{\Delta}=H_{\Gamma}\times \R$ and $U_{\Delta}=U_{\Gamma}\times\R$.
As above we have the quasi-isomorphism inclusion $ (\bigwedge {\mathfrak u} ^{\ast})^{\{\pm 1\}}\otimes \bigwedge (y)\subset A^{\ast}(M_{\Delta})$.
Let $\omega=x_{1}\wedge y+x_{2}\wedge x_{3}$.
Then $\omega$ is a symplectic form on $M_{\Delta}$.
Since $H^{1}(M_{\Delta},\R)\cong \langle x_{1}, y \rangle$ and $H^{3}(M_{\Delta},\R)\cong \langle x_{1}\wedge x_{2} \wedge x_{3}, x_{2}\wedge x_{3}\wedge y\rangle $, the linear map $[\omega]\wedge : H^{1}(M_{\Delta},\R)\to H^{3}(M_{\Delta},\R)$ is an isomorphism and hence we have the following proposition.
\begin{proposition}
$M_{\Gamma}\times S^{1}$ satisfies the hard Lefschetz property.
\end{proposition}
\begin{remark}
$\Delta$ is a finite extension group of the non-abelian nilpotent group $2\Z\ltimes\Z^{2}\times \Z$ as remark \ref{fi}.
By the result of \cite{BG}, a compact $K(2\Z\ltimes\Z^{2}\times \Z, 1)$-manifold  is not a Lefschetz $4$-manifold.
Thus $M_{\Delta}$ is a example of a Lefschetz $4$-manifold with non-Lefschetz finite covering space.
In \cite[Example 3.4]{Li}, Lin showed  the existence of Lefschetz $4$-manifolds with non-Lefschetz finite covering space.
$M_{\Delta}$ is a simpler and more constructive example.
\end{remark}

\ \\
{  Acknowledgements.} 

The author would like to express his gratitude to   Toshitake Kohno for helpful suggestions and stimulating discussions.
He would also like to thank     Katsuhiko Kuribayashi and Keizo Hasegawa for their active interests in this paper.
This research is supported by JSPS Research Fellowships for Young Scientists.


\begin{thebibliography}{20}
\bibitem{A}D. Arapura,  K$\ddot {\rm a}$hler solvmanifolds, Int. Math. Res. Not.{\bf 3}(2004), 131-137.
\bibitem{Aus}
L. Auslander, An exposition of the structure of solvmanifolds. I. Algebraic theory.  Bull. Amer. Math. Soc. {\bf 79}  (1973), no. 2, 227--261.

\bibitem{B}
O. Baues, Infra-solvmanifolds and rigidity of subgroups in solvable linear algebraic groups. Topology {\bf 43} (2004), no. 4, 903--924.
\bibitem{BC} O. Baues,  V. Cort\'es,
Aspherical K\"ahler manifolds with solvable fundamental group. Geom. Dedicata {\bf 122} (2006), 215--229.


 \bibitem{BG}
C. Benson, and C. S. Gordon, 
K$\ddot {\rm a}$hler and symplectic structures on nilmanifolds. 
Topology {\bf 27} (1988), no. 4, 513--518. 
\bibitem{Boc}
C. Bock, On Low-Dimensional Solvmanifolds, preprint arXiv:0903.2926 (2009).

\bibitem{Bor}
A. Borel, Linear algebraic groups 2nd enl. ed Springer-verlag (1991).
\bibitem{Bur} D. Burde, K. Dekimpe, K. Vercammen,
Complete LR-structures on solvable Lie algebras. J. Group Theory 13 (2010), no. 5, 703--719.

\bibitem{DGMS}
P. Deligne, P. Griffiths, J. Morgan, and D. Sullivan, 
Real homotopy theory of Kahler manifolds. Invent. Math. {\bf 29} (1975), no. 3, 245--274.
\bibitem{FG}
M. Fernandez, and A. Gray,  Compact symplectic solvmanifolds not admitting complex structures. Geom. Dedicata {\bf 34} (1990), no. 3, 295--299. 
	\bibitem{H} K. Hasegawa, Minimal models of nilmanifolds. Proc. Amer. Math. Soc. {\bf 106} (1989), no. 1, 65--71. 

\bibitem{Hn} K. Hasegawa, A note on compact solvmanifolds with K\"ahler structures. Osaka J. Math. {\bf 43} (2006), no. 1, 131--135.
\bibitem{Hd} K. Hasegawa, Small deformations and non-left-invariant complex structures on six-dimensional compact solvmanifolds. Differential Geom. Appl. {\bf 28} (2010), no. 2, 220--227.
\bibitem{Hatt} A. Hattori, Spectral sequence in the de Rham cohomology of fibre bundles. J. Fac. Sci. Univ. Tokyo Sect. I {\bf 8} 1960 289--331 (1960). 
\bibitem{Na} I. Nakamura, Complex parallelisable manifolds and their small deformations. J. Differential Geometry {\bf 10} (1975), 85--112.
\bibitem{Li} Y. Lin, Examples of non-Kahler Hamiltonian circle manifolds with the strong Lefschetz property. Adv. Math. {\bf 208} (2007), no. 2, 699--709.
\bibitem{OT}K. Oeljeklaus, M. Toma,  Non-K\"ahler compact complex manifolds associated to number fields. Ann. Inst. Fourier (Grenoble) {\bf 55} (2005), no. 1, 161--171.
\bibitem{OV2}
A. L. Onishchik, and E. B. Vinberg, Structure of Lie groups and Lie algebras, Springer (1994).
\bibitem{OV}
A. L. Onishchik, and E. B. Vinberg, Discrete subgroups of Lie groups and cohomologies of Lie groups and Lie algebras, Springer (2000).
 \bibitem{T}
J. Oprea, and A. Tralle, Symplectic manifolds with no K$\ddot {\rm a}$hler structure. Lecture Notes in Math. 1661, Springer (1997).

 \bibitem{R}
 M. S. Raghnathan, Discrete subgroups of Lie Groups, Springer-verlag, New York, 1972.
 \bibitem{Re}
 B. E. Reed,  Representations of solvable Lie algebras, Michigan Math. J.  {\bf 16} 1969 227--233. 
 
 \bibitem{PS} 
 P. Scott,  The geometries of $3$-manifolds. Bull. London Math. Soc. {\bf 15} (1983), no. 5, 401--487

\end{thebibliography}
\end{document}